\documentclass[a4paper,12pt,reqno]{amsart}
\usepackage{amssymb,delarray,amsmath,stmaryrd}
\usepackage[dvips]{graphicx}

\usepackage{hyperref}
\usepackage[all]{xy}

\input epsf.tex
\newdimen\xsize
\newdimen\oldbaselineskip
\newdimen\oldlineskiplimit
\def\restorelineskip{\baselineskip=\oldbaselineskip%
\lineskiplimit=\oldlineskiplimit}
\def\putm[#1][#2]#3{
\hbox{\vbox to 0pt{\parindent=0pt%
\vskip#2\xsize\hbox to0pt{\hskip#1\xsize $#3$\hss}\vss}}}%
\long\def\Line#1{\hbox to \hsize{#1}}
\def\putt[#1][#2]#3{
\vbox to 0pt{\noindent\hskip#1\xsize\lower#2\xsize%
\vtop{\restorelineskip#3}\vss}}

\makeatletter
\def\xbig[#1]#2{{\hbox{$\m@th\left#2\vbox to#1\xsize{}%
\right.\n@space$}}}
\makeatother
\def\xlar[#1]#2{%
\smash{\mathop{ \hbox to #1\xsize{\leftarrowfill}}\limits^{#2}}}
\def\xrar[#1]#2{%
\smash{\mathop{ \hbox to #1\xsize{\rightarrowfill}}\limits^{#2}}}
\def\xline[#1]{\hbox to #1\xsize{\leaders\hrule\hfill}}

\DeclareFontFamily{U}{rsf}{\skewchar\font'177}%
\DeclareFontShape{U}{rsf}{m}{n}{<-6>rsfs5<6-8>rsfs7<8->rsfs10}{}%
\DeclareFontShape{U}{rsf}{b}{n}{<-6>rsfs5<6-8>rsfs7<8->rsfs10}{}%
\DeclareMathAlphabet\RSFS{U}{rsf}{m}{n}
\SetMathAlphabet\RSFS{bold}{U}{rsf}{b}{n}
\def\sf#1{{\mathsf{#1}}}

\def\slsf{\slshape \sffamily }

\def\msmall#1{\mathchoice{\hbox{\small$\displaystyle {#1}$}}{#1}{#1}{#1}}

\def\cc{{\mathbb C}}
\def\dd{{\mathbb D}}

\def\rr{{\mathbb R}}

\def\st{_{\mathsf{st}}}

\def\id{\sf{Id}}

\def\lim{\mathop{\sf{lim}}}

\def\log{\sf{log}\,}

\def\eps{\varepsilon}

\def\<{\langle}\let\la=\<
\def\>{\rangle}\let\ra=\>
  
\def\comp{\Subset}
\def\d{\partial}

\def\ddef{\mathrel{{=}\raise0.3pt\hbox{:}}}
\def\deff{\mathrel{\raise0.3pt\hbox{\rm:}{=}}}

\def\fraction#1/#2{\mathchoice{{\msmall{ #1\over#2}}}%
{{ #1\over #2 }}{{#1/#2}}{{#1/#2}}}
\def\norm#1{\left\Vert{#1}\right\Vert}

\def\le{\leqslant}

\def\emptyset{\varnothing}

\def\longpoints{\leaders\hbox to 0.5em{\hss.\hss}\hfill \hskip0pt}
\def\stateskip{\smallskip}
\def\state#1. {\stateskip\noindent{\bf#1. }} 
\def\statep#1. {\stateskip\noindent{\bf#1 }} 

\def\Chi{\raise 2pt\hbox{$\chi$}}

\def\ie{\hskip1pt plus1pt{\sl i.e.\/,\ \hskip1pt plus1pt}}

\def\Chi{\raise 2pt\hbox{$\chi$}}

\let\phI=\phi\let\phi=\varphi\let\varphi=\phI
\let\cal=\mathcal

\def\calb{{\cal B}}
\def\calc{{\cal C}}

\def\eps{\varepsilon}

\def\comp{\Subset}
\def\d{\partial}

\def\1{{1\mkern-5mu{\rom l}}}

\def\fraction#1/#2{\mathchoice{{\msmall{ #1\over#2}}}%
{{ #1\over #2 }}{{#1/#2}}{{#1/#2}}}

\def\le{\leqslant}

\def\emptyset{\varnothing}

\def\qed{\ \ \hfill\hbox to .1pt{}\hfill\hbox to .1pt{}\hfill $\square$\par}

\def\comment#1\endcomment{}

 \belowdisplayskip=\abovedisplayskip

\def\lineeqqno(#1){\hfill\llap{\vbox to 10pt%
{\vss\begin{align} \eqqno(#1)\end{align}\vss}}\vskip1pt}

\advance\headheight 1.2pt

\def\ShowwLLabel#1{}

\def\thechpt{\Roman{chpt}}

\def\newchapt[#1]#2{\newpage%
\refstepcounter{chpt}\setcounter{subsection}{0}%
\setcounter{thm}{0}\setcounter{defi}{0}%
\setcounter{rema}{0}\setcounter{exrc}{0}%
\renewcommand{\thesubsection}{\thechpt.\arabic{subsection}}%
\section*{\begin{center}\huge \bf Chapter \thechpt\\
#2 \end{center}}\label{#1}%
\ \smallskip%
\markboth{Chapter \thechpt}{#2}%
}
\def\newsect[#1]#2{\refstepcounter{section}\setcounter{equation}{0}%
\renewcommand{\thesubsection}{\arabic{section}.\arabic{subsection}}%
\section*{\arabic{section}.
#2}\vspace{-20pt}\label{#1}\vspace{20pt}%
\markboth{Section \arabic{section}}{#2}}

\def\newlect[#1]#2{\refstepcounter{section}%
\renewcommand{\thesubsection}{\arabic{section}.\arabic{subsection}}%
\section*{Lecture \arabic{section}\\
#2}\label{#1}%
\markboth{Lecture \arabic{section}}{#2}}

\def\newprg[#1]#2{\refstepcounter{subsection}%
\subsection*{{\thesubsection.\ #2}} \label{#1}%
}

\def\newappx[#1]#2{%
\refstepcounter{appx}\setcounter{section}{0}%
\renewcommand{\thesubsection}{A\arabic{appx}.\arabic{subsection}}%
\section*{Appendix \arabic{appx}\\ #2}
\label{#1}%
\markboth{Appendix A\arabic{appx}}{#2}
}

\newtheorem{thm}{Theorem}[section]
   \def\newthm#1{\begin{thm}\label{#1}}

\newtheorem{nnthm}{Theorem.} 
   \def\newnnthm#1{\begin{nnthm} \label{#1}}
\newtheorem{lem}{Lemma}[section]
   \def\newlemma#1{\begin{lem} \label{#1}}

\newtheorem{prop}{Proposition}[section]
   \def\newprop#1{\begin{prop}\label{#1}}

\newtheorem{corol}{Corollary}[section]
   \def\newcorol#1{\begin{corol} \label{#1}}

\newtheorem{defi}{Definition}[section]
   \def\newdefi#1{\begin{defi} \label{#1}\rm }

\newtheorem{exmp}{Example}[section]
   \def\newexmp#1{\begin{exmp} \label{#1}\rm }

\newtheorem{exrc}{Exercise}
   \def\newexrc#1{\begin{exrc} \label{#1}\rm }

\newtheorem{quest}{Open Question}
   \def\newrema#1{\begin{quest} \label{#1}\rm }

\newtheorem{rema}{Remark}
   \def\newrema#1{\begin{rema} \label{#1}\rm }

\let\xrar=\xrightarrow

\def\eqqno(#1){\label{(#1)}}
\def\eqqref(#1){(\ref{(#1)})}

\numberwithin{equation}{section}

\title{Boundary Values and Boundary Uniqueness of $J$-Holomorphic Mappings}
\author{S. Ivashkovich, J.-P. Rosay}
\address{
Universit\'e de Lille-1, UFR de Math\'ematiques, 59655 Villeneuve
d'Ascq, France.}
\email{ivachkov@math.univ-lille1.fr}
\address{IAPMM Nat. Acad. Sci. Ukraine,
Lviv, Naukova 3b,
79601 Ukraine.}
\address{
Department of Mathematics, University of Wisconsin, Madison WI 53706 USA. }
\email{ jrosay@math.wisc.edu}

\date{\today}
\begin{document}

\begin{abstract}
We establish a Fatou-type Theorem for $J$-holomorphic mappings that
are bounded in an appropriate sense and we prove the Blaschke condition
for their zero sets. We also prove a Privalov-type uniqueness Theorem for
pairs of $J$-holomorphic mappings.

\end{abstract}
\maketitle

\tableofcontents

\smallskip\qed
\section[1]{Introduction}

Bounded holomorphic functions, defined in the unit disc
$\Delta\subset \cc$, have special properties. Namely, if $f$ is such
a function then the radial limit $f^*(e^{i\theta})
\deff\lim_{r\to 1}f(re^{i\theta})$  exists almost everywhere and, if
$f\not\equiv 0$, then the set of  $\theta\in [0, 2\pi ]$ such that $f^\ast(e^{i\theta}) = 0$ has zero length.
Further, the zero set $\{\zeta_k\}\subset \Delta $ of a bounded holomorphic function
satisfies the Blaschke condition.

\smallskip
Our goal in this paper is to generalize these results to the case of
pseudoholomorphic mappings with values in almost complex manifolds. The first problem that one meets here
is to introduce a correct notion of boundedness. We propose the following one:
\begin{itemize}
\item A $J$-holomorphic mapping $u:\Delta\to X$ into an almost complex manifold
$(X,J)$ is called {\slsf bounded} if there exists a relatively
compact domain  $\Omega\supset\overline{u(\Delta)}$ which admits a
strictly $J$-plurisubharmonic function.
\end{itemize}
Note that this function can be always supposed to be bounded from above. In the classical function
theory $|\zeta^2|$ plays the role of this function.

\smallskip Likewise, a relatively compact domain  $\Omega\comp X$ is called {\slsf bounded}
if it admits a strictly $J$-plurisubharmonic function bounded from
above. The definition of boundedness of holomorphic maps or domains,
which we propose, requires the existence of a bounded strictly
$J$-plurisubharmonic function. Neither asking relative compactness
nor the existence of a non-constant bounded $J$-plurisubharmonic
function would suffice. In both cases, there could be $J$-complex
lines, i.e. non-constant images of $\cc$ under a $J$-holomorphic
mapping. But that should clearly be forbidden.

\smallskip We shall see that bounded $J$-holomorphic mappings have
non-tangential limits for almost all $e^{i\theta}\in\d\Delta$, \ie
that a Fatou-type theorem is valid for them.

\begin{thm}
\label{fatou} Let $u$ be a bounded $J$-holomorphic map from $\Delta$
into an almost complex manifold $(X,J)$, $J$ continuous.
Then for almost every $e^{i\theta} \in \d\Delta$ the limit
$u^*(e^{i\theta})$ of $u(\zeta )$ exists as $\zeta$ approaches
$e^{i\theta}$ non-tangentially.
\end{thm}

\medskip Our second result is the following.

\begin{thm}
\label{bla}
Let $u:\Delta\to X$ be a non-constant, bounded $J$-holomorphic
mapping, $J\in\calc^{Lip}$. Then for every point $p\in X$ the
set $\{\zeta_k\in \Delta : u(\zeta_k) = p\}$ satisfies the Blaschke
condition
\begin{equation}
\label{blasch}
\sum_k(1-|\zeta_k|) < +\infty .
\end{equation}
\end{thm}
Here we denote by  $\calc^{Lip}$ the class of Lipschitz-continuous functions, and by
$\calc^{1,Lip}$ the class of functions whose first derivatives are Lipschitz-continuous.

\smallskip  We prove in this paper the following direct
analogue of the Privalov boundary uniqueness theorem.

\begin{thm}
\label{priva} Suppose $J\in\calc^{1,Lip}$ and let $u_1,u_2:\Delta\to
X$ be $J$-holomorphic mappings. Assume that there is a set $E$ of
positive measure on the unit circle, such that at every $
e^{i\theta}\in E$ $u_1$ and $u_2$ have non-tangential boundary
values and that these boundary values are equal, \ie $u_1^*(
e^{i\theta}) = u_2^*( e^{i\theta})$ for $e^{i\theta}\in E$. Then
$u_1\equiv u_2$.
\end{thm}
The result is still valid with a weaker notion of non-tangential
boundary value (called restricted in \S 3), see Theorem
\ref{geo-incl} below.

\medskip\noindent {\slsf Notes.} 1) We would like to draw the
attention of the reader to the fact that we do not require the
continuity of the strictly $J$-plurisubharmonic function in the definition of
boundedness of functions and domains. Assuming continuity would
allow a simplification in the proof of Theorem \ref{fatou}. However,
it is not at all clear that the existence of non-continuous strictly
$J$-plurisubharmonic function insures the existence of continuous
ones, nor is it clear that the existence of continuous ones insures
the existence of smooth ones. These are challenging problems.

\smallskip\noindent 2) Let us emphasize that the statement of Theorem
\ref{bla} is not as satisfactory as the statement of Theorem \ref{priva}.
One should try to get the Blaschke condition (\ref{blasch})
for the set of coincidence $\{\zeta_k : u_1(\zeta_k)=u_2(\zeta_k)\}$
of two distinct bounded $J$-holomorphic mappings. However, for the
reasons explained in Section 4, this general statement is out of the
range of the methods deployed in this paper which exploits
plurisubharmonic functions and the notion of {\slsf pluripolarity}.
That is the main tool that allows us to go beyond the preliminary
uniqueness results established in \cite{ISk}.

\smallskip\noindent 3) Boundary uniqueness problems
for solutions of certain PDE-s (mostly close to the Laplace
equation) are commonly studied with an {\slsf ad hoc} assumption of
some (mostly $\calc^1$) regularity up to the boundary, see \cite{AB}
for example. It is important to point out that we do not assume
{\slsf any}  boundary regularity of our mappings and in general there is
no regularity.

\section[2]{Plurisubharmonic Functions and non-Tangential Limits}

\smallskip\noindent{\slsf 2.1 Plurisubharmonic functions}

\smallskip An upper semi-continous function $\rho$ defined
on an almost complex manifold $(X,J)$, with $J$ continuous,
is said to be $J$-plurisubharmonic ($J$-psh for short) if its restriction to any
$J$-holomorphic disc is subharmonic; i.e. if $u:\Delta \to (X,J)$
is $J$-holomorphic, then $\rho\circ u$ is subharmonic. Of course this defintion is
meaningful in case there is abundance of $J$-holomorphic discs (e.g. the maximum principle certainly holds),
and this happen as soon as $J$ is {\slsf continuous}, see Appendix.

\smallskip
We say that $\rho$ is strictly $J$-psh if and only if, locally, $\rho+\eta$ is still $J$-psh for all functions $\eta$ of small
enough $\calc^2$ norm.

\smallskip
If the almost complex structure  $J$ is at least of class $\calc^1$,
there is for $\calc^2$ functions $\rho$ a characterization of $J$-plurisubharmonicity in terms of the
Hessian (Levi form). The condition is:
\[
dd^c_J\rho_p(v,J(p)v)\geq 0
\]
for any tangent vector $v$ at any point $p\in X$. Here, for a function $\varphi$,
$d^c\varphi (v)=-d\varphi(Jv)$ for $v\in TX$.
This is rather staightforward if $J$ is of class at least $\calc^{1,\alpha}$ for some $0<\alpha<1$.
The $J$-holomorphic discs are then of class $\calc^2$ (and even $\calc^{2,\alpha}$). More care is needed
in case of $J$ merely $\calc^1$ since the $J$-holomorphic discs may fail to be of class $\calc^2$.
See Corollary 1.1 in \cite{IR}. For strict $J$-plurisubhamonicity the condition is that
$dd^c_J\rho(v,Jv)>0$ if $v\neq 0$. Equivalently the conditions can be given in terms of
$\partial_J\overline\partial_J$ as in \cite{P}.

\smallskip\noindent{\slsf 2.2. Non-tangential limits}

\smallskip We start with the study of  non-tangential boundary values of a bounded
$J$-holomorphic map from the unit disc $\Delta $ in $\cc$ with
values in an almost complex manifold $(X,J)$, with $J$ merely continuous.
We shall consider domains in the unit disc obtained by
taking conic neighborhoods of radii. Adopting Rudin's notation from \cite{Ru} (11.18),
for $0< \alpha < 1$ and $\theta\in [0, 2\pi]$, we consider the set
\begin{equation}
\eqqno(conic)
e^{i\theta}\Omega_\alpha = \{\zeta\in\Delta : |\zeta - e^{i\theta}|\zeta || < \alpha (1 - |\zeta |)\}.
\end{equation}
This set is a conic neighborhood of the radius $[0,e^{i\theta}]$ with vertex at $e^{i\theta}$. Let
$u:\Delta\to X$ be a $J$-holomorphic map.
\begin{defi}
\label{ntg-lim} We say that $u$ has a non-tangential limit at
$e^{i\theta}$ if for every $\alpha\in (0,1)$ $u(\zeta )$ has a limit
as $\zeta\to e^{i\theta}$, $\zeta \in e^{i\theta}\Omega_\alpha$.
This limit is then denoted as $u^\ast (e^{i\theta})$.
\end{defi}

\medskip\noindent{\sl Proof of the Theorem \ref{fatou}.}
We start the proof with two preliminaries, one on some kind of
Schwarz Lemma, and one on radial limits a.e. for subharmonic
functions. Let $\Omega$ be a relatively compact domain in $X$
containing $\overline{u(\Delta)}$ and $\rho$ a strictly $J$-psh
function in $\Omega$. Fix some Riemannian metric $h$ on $X$. We
claim that for any compact set $K\comp \Omega$ there exists a
constant $C =C(K,h,\rho)$ such that for any $J$-holomorphic
$v:\Delta\to K$ one has
\begin{equation}
\label{schwarz} \norm{dv(0)}_h\le C.
\end{equation}
Otherwise, one gets a sequence of $J$-holomorphic $v_n:\Delta\to K$
with $\norm{v_n(0)}_h\to\infty$. Then, the Brody re-parameterization
Lemma (in which neither the integrability of $J$ nor its smoothness play any
role) \cite{B} applies and gives a non-constant $J$-holomorphic map $v:\cc\to
K$. The function $\rho\circ v$ is then a bounded subharmonic
function, hence it is constant. Since $\rho$ is strictly
plurisubharmonic this is impossible.

\smallskip
Remark that by simple rescaling, (\ref{schwarz}) gives
$\norm{dv(\zeta )}_h\le \frac{C}{1-|\zeta |}$, for $v$ as above, and
it follows that
\begin{equation}
\label{schwarz2}
{\rm dist}_h(v(\zeta ),v(\zeta ')\leq C
\delta (\zeta ,\zeta'),
\end{equation}
where $\delta$ denotes the Poincar\'e metric on the unit disc.

\smallskip The second essential tool in the proof of Theorem \ref{fatou} is
a Theorem of Littlewood \cite{L} that asserts that any subharmonic
function on the unit disc that is bounded form above, has radial
boundary values a.e., see also \cite{T}, IV.10. It is
extremely easy to see that in general there is no non-tangential
limits a.e. Under hypotheses of normality non-tangential boundary
values however exist, see \cite{M}. If we assumed continuity of
$\rho$, \cite{M} would apply due to (\ref{schwarz2}). We will simply
have to use the argument of normality, one step later.

\medskip
First, we state a generalization of Littlewood's Theorem (where
instead of taking limits along radii one takes limits along rays). A
very convenient reference (where much more is done) is \cite{D}. Let
$\lambda$ be a subharmonic function defined on $\Delta$ that is bounded
from above. Let $\nu\in (-\frac{\pi }{2}, \frac{\pi }{2})$, we shall
say that $\lambda$ has boundary at $e^{i\theta}$ in the $\nu$-direction
if $\lambda (e^{i\theta}-te^{i(\theta -\nu)})$ has a limit as
$t\searrow 0$. Denote this limit by $\lambda^\ast_{\nu}(e^{i\theta})$.
Of course $\nu =0$ corresponds to the radial limit. The result of
\cite{D} is that for any fixed $\nu$ the limit
$\lambda^\ast_{\nu}(e^{i\theta})$ exists a.e. Moreover if $\nu'\neq
\nu$, $\lambda^\ast_{\nu}=\lambda^\ast_{\nu'}$ a.e. The last fact is
because the heart of the matter is to show that a Green potential
has boundary value 0 a.e. (in the direction of $\nu$). See Corollary
1 on p. 520 in \cite{D}.

\smallskip
After these preliminaries we are ready for the proof.
Let $\rho$ be a strictly $J$-psh function defined on a neighborhood of
$\overline{u(\Delta )}$. Taking $\lambda=\rho\circ u$ above, we see that
if $\nu\in
(-\frac{\pi }{2}, \frac{\pi }{2})$, $\rho\circ u$ has boundary value in
the $\nu$-direction a.e.

\smallskip
Let $\phi =(\phi_1,\cdots ,\phi_N)$ be a $\calc^1$ injective map
from a neighborhood of $\bar\Omega$ into $\rr^N$ (with $N$ large
enough). We can apply the same reasoning to the functions
$\lambda_{\eps}\deff\rho\circ u + \epsilon \phi_j\circ u$,
$j=1,...,N$ with $\eps >0$ small enough. It follows that $u$ has
boundary value $u^\ast_\nu$ a.e in the $\nu$-direction, and as
discussed above the boundary values for any two values of $\nu$ will
agree a.e.

\smallskip
In order to get non-tangential boundary values for $u$ we have to
use a normality argument now. Let $F$ be a dense countable subset of
$(-\frac{\pi }{2}, \frac{\pi }{2})$. There exists a set $\Theta$ of full
measure in $[0,2\pi)$ such that for all $\nu\in F$ and all $\theta
\in \Theta$, $u^\ast_\nu(e^{i\theta})$ exists and does not depend on
$\nu$, we thus simply write $u^\ast(e^{i\theta})$. At any such point
$e^{i\theta}$, $u$ has a non-tangential boundary value. Indeed, if
$0< \alpha <1$ and $\eps >0$ are fixed, there are $\nu_1,\cdots
,\nu_N\in F $ such that the the union of
the rays from $e^{i\theta}$ in the $\nu_j$-directions ($j=\{
1,\cdots ,N\}$) is $\epsilon$ dense in the Poincar\'e metric in
$e^{i\theta}\Omega_\alpha$. There exists $\eta >0$ such that
for $\zeta$ on any of these $N$ rays in the $\nu_j$-directions
${\rm dist}~(u(\zeta),u^\ast
(e^{i\theta}))\leq \epsilon$, if $1-\eta<|\zeta|<1$.
Then (\ref{schwarz2}) shows that for
$\zeta\in e^{i\theta}\Omega_\alpha$, ${\rm dist}~(u(\zeta),u^\ast
(e^{i\theta}))\leq (C+1) \epsilon$, if $1-\eta <|\zeta|$.

\smallskip\qed

\smallskip\noindent{\slsf 2.3. Pluripolarity}

\smallskip\noindent
Pluripolarity refers to the $-\infty$ sets of $J$-psh functions. More precisely,
a subset $V$ of an almost complex manifold $(X,J)$, $J$ continuous, is called {\slsf
locally complete pluripolar} if for any point $p\in X$ there exists
a neighborhood $U\ni p$ and a $J$-psh function $\rho$ in $U$, not
identically equal $-\infty$, such that $M\cap U = \{ x\in U: \rho
(x) = -\infty\}$. If as $U$ can be taken a neighborhood of $V$ then
$V$ is called {\slsf complete pluripolar}, or globally complete
pluripolar.

\smallskip The following examples of pluripolarity require
some smoothness of $J$.

\medskip\noindent
- Chirka's function: If $J$ is a $\calc^1$-regular almost
complex structure defined near 0 in $\rr^{2n}=\cc^n$, and if $J(0)=J_{st}$
(the standard complex structure), then for $A>0$, large enough, the
function $\log |z|~+~A|z|$, is $J$-psh near 0 (\cite{IR} Lemma 1.4). That implies that a
point is a complete pluripolar set. The result extends to almost complex structures that are
Lipschitz-continuous, see \S 5.2   below.

\medskip\noindent
- Elkhadhra's theorem: Let $J$ be a $\calc^{1,Lip}$-regular almost
complex structure defined near 0 in $\rr^{2n}$ and let $V$ be a germ
at zero  of a $J$-complex submanifold of class $\calc^{2,Lip}$. Then
there exist a neighborhood of $U$ of 0 and a $J$-psh function $\rho$
defined in $U$ such that $\rho$ is of class $\calc^{2,Lip}$ on
$U\setminus V$ and $V\cap U = \rho^{-1}(-\infty )$, see \cite{E}.
The function  $\rho$ is of the form $-\log \big(-\log |f|^2\big) + A
|z|^2 $, where $z=(z_1,...,z_n)$ are local complex coordinates
centered at $0$ such that $J(0)=J\st$ in these coordinates  and
$|f|^2=\sum_{j=1}^p|f_j|^2$, for appropriate functions $f_j$'s
vanishing on $V$. I.e., $V$ is a "$\log\log$"-polar set of a $J$-psh
function, unlike the case of a point, where the Chirka function
has a "$\log$"-pole. Elkhadhra's Theorem generalizes the result of
\cite{Ro1} for $J$-holomorphic curves in which already a
``$\log\log$'' function was used.
The difference between $\log$ singularities and lesser singularities
(that is essentail in classical complex pluri-potential theory)
plays an important role in this paper, see Section 4.

\smallskip\noindent{\slsf Note .} The proof of Elkhadhra's theorem
requires
$\calc^{1,Lip}$ smoothness of $J$.
In \cite{Ro2} (where a correction for the smoothness
hypothesis in \cite{Ro1} is made) a possibly more transparent proof
of pluripolarity of curves is given. That proof can be adapted to
give a different proof of Elkhadhra's Theorem, but it also very
clearly requires $\calc^{1,Lip}$ regularity of $J$. It is clear that
$J$-holomorphic curves can fail to be locally complete pluripolar if the almost complex
structure is only of class $\calc^\alpha$ (for some $0<\alpha <1$). Indeed two
$J$-holomorphic maps $u_1$ and $u_2$, both imbeddings, can be such that they do not coincide
on any neighborhood of 0, but they coincide on an open set whose closure contains 0.
Any $J$-psh function $\rho$  with value $-\infty$ on one of the discs near $u_1(0)=u_2(0)$
will also be $-\infty$ on the other disc. So the case left unclear is the case of
structures $J\in\calc^{Lip}$ (or, even of smoothness $\calc^1$ or $\calc^{1,\alpha}$, $0<\alpha <1$).
Are $J$-holomorphic curves still locally complete pluripolar and more generally does Elkhadhra's theorem extend
when $J$ is only Lipschitz-continuous?

\section[3]{Geometric Inclusion and Privalov-type Theorem}

\smallskip\noindent{\slsf 3.1. Stratified pluripolar sets}

\medskip A Privalov-type uniqueness theorem will be obtained in this paper as an immediate corollary
of a more general and more "geometric" statement, we call it the
"geometric inclusion". Let us start with the following

\begin{defi}
\label{strat} Let $(X,J)$ be an almost complex manifold, $J$ continuous. Let
$V\subset X$ be a closed subset. We say that $V$ admits a {\slsf
locally complete pluripolar stratification} if $V$ can be
represented as the union $V=V_{k}\cup V_{k-i}\cup ... \cup V_{0}$
where $V_{0}$ is closed, locally complete pluripolar in $(X,J)$ and
$V_{j}$ is closed, locally complete pluripolar in $X\setminus
\bigcup_{l=0}^{j-1}V_{l}$ for $j=1,...,k$.
\end{defi}
In order to justify this definition let us give two examples.

\smallskip\noindent{\bf 1.} Any $\calc^{2,Lip}$ $J$-complex submanifold in $(X,J)$ is locally complete pluripolar,
provided that $J\in \calc^{1,Lip}$, see above. Especially we shall use the fact that the diagonal $\dd$
in the Cartesian square
$(X\times X, J\oplus J)$ of an almost complex manifold is locally complete pluripolar, if $J$ is of class
$\calc^{1,Lip}$.

\smallskip\noindent{\bf 2.} If $J$ is an almost complex structure
of class $\calc^{2,\alpha}$ and $V$ is a closed $J$-complex curve in
$X$ then the singular part $V_0$ of $V$ is discrete and therefore
locally complete pluripolar. $V_1\deff V\setminus V_0$ is of class $\calc^{3,\alpha}$
(we only need $\calc^{2,Lip}$) and thus
is  locally the $-\infty$ set of a $J$-psh function.
Therefore  a $J$-complex curve admits a
locally complete pluripolar stratification.

\smallskip\noindent{\bf 3.} These two examples could be generalized as follows.
Here for simplicity we shall just assume the data to be smooth.
Let us define a
$J$-analytic subset of an almost complex manifold $(X,J)$ as a
stratified set $A=A_k\cup A_{k-1}\cup ...\cup A_0$, where $A_0$ is a
$J$-complex submanifold in $X$ and $A_j$ is a $J$-complex
submanifold of $X\setminus \bigcup_{l=0}^{j-1}A_l$. A motivation for
such a definition obviously comes from the structure theorem for the
usual analytic sets. From Elkhadra's Theorem it follows that a
$J$-analytic set is locally complete pluripolar stratified.

\medskip\noindent{\slsf 3.2. Geometric inclusion}

\medskip\noindent
In \S 2.1 we have given boundedness hypotheses that guarantee the existence of
non-tangential limits. In this section, we proceed differently. We shall assume the existence
of non-tangential limits but we will no longer require that our
mappings are bounded. Moreover, unlike in  Definition \ref{ntg-lim}
we shall need only the existence of $\lim u(\zeta )$, as $\zeta$
approaches to $e^{i\theta}$, $\zeta\in e^{i\theta}\Omega_\alpha$,
for some $\alpha$ possibly depending on $\theta$. In that case we say
that $u$ has a restricted non-tangential limit at $e^{i\theta}$.

\begin{thm}
\label{geo-incl}
Let $V$ be a locally complete pluripolar stratified subset of an almost complex manifold $(X,J)$, with
$J$ continuous, and
let $u:\Delta\to X$ be a $J$-holomorphic map. Assume that for $e^{i\theta }$ in a set $E$ of positive
measure in $\d\Delta$ the mapping $u$ has a restricted non-tangential limit $u^*(e^{i\theta})$ and that
$u^*(e^{i\theta})\in V$. Then $u(\Delta)\subset V$.
\end{thm}

\noindent Proof. For a set $E\subset \d\Delta$ and $0<r<1$, let
\begin{equation}
\eqqno(gamma)
\Gamma_r(E,\alpha )= \bigcup_{e^{i\theta}\in E} e^{i\theta}\Omega_\alpha~~\cap \{ |\zeta |>r\}~.
\end{equation}
This set is therefore made of truncated cones of fixed aperture with
vertices at the points $e^{i\theta}\in E$.  In the sequel we fix the
minimal $j$ such that the set of $e^{i\theta}$ with
$u^*(e^{i\theta})\in V_j\setminus \bigcup_{l=0}^{j-1}V_l$ is of
positive length. We denote this set still by $E$. If $j=0$ we denote
by $E$ the set of positive length such that $u^*(e^{i\theta})\in
V_0$ for $e^{i\theta}\in E$.

\medskip\noindent{\slsf Step 1.} {\sl Let $E$ be as just defined. Then there exist
$p\in V_{j}\setminus\bigcup_{l=0}^{j-1}V_l$, a neighborhood $W$ of
$p$ in $X\setminus \bigcup_{l=0}^{j-1}V_l$ and a $J$-psh function
$\rho$ defined on $W$, bounded from above, with $V_{j}\cap
W=\rho^{-1}(-\infty )$, and there exist a closed set $E_0\subset
\d\Delta$ of positive length, $r\in (0,1)$, and $\alpha\in (0,1)$
such that $u(\Gamma_r(E_0,\alpha ))$ is a relatively compact subset
of $W$, and for all $e^{i\theta}\in E_0$ $u^\ast (e^{i\theta})\in
V_{j}\cap W$.}

\smallskip\noindent The mapping $u^*:E\to X$ is measurable as a pointwise limit of measurable mappings.
By Lusin's theorem,  we can assume that $E$ is a closed set and that
$u^\ast$ is continuous on $E$. Fix $e^{i\theta_0}\in E$ such that
$E$ has positive length in any neighborhood of $e^{i\theta_0}$, set
$p=u^\ast (e^{i\theta_0})$. Let $W_0$ and $W$ be neighborhoods of
$p$, $W_0\subset\subset W$, $W\cap \bigcup_{l=0}^{j-1}V_l=\emptyset$
such that there exists a $J$-psh function $\rho$ defined on $W$ such
that $V_j\cap W$ coincides with $\rho^{-1}(-\infty )$ in $W$.
Restricting $E$ to a small neighborhood of $e^{i\theta_0}$ we can
assume that for $e^{i\theta}\in E$ one has $u^\ast(e^{i\theta})\in
W_0$. For each $e^{i\theta}\in E$ there exists $\alpha>0$ and $r<1$
such that $u( e^{i\theta}\Omega_\alpha~~\cap \{ |\zeta |>r\})\subset
W_0\subset\overline{W_0}$, where $\alpha$ and $r$ depend on
$\theta$. To finish the proof of the claim we simply need to have
$\alpha$ and $r$ no longer depending on $\theta$. For any integer
$k>0$, set
\[
E_k=\{ e^{i\theta}\in E~: u(e^{i\theta}\Omega_{\frac{1}{k}}\cap \{ |\zeta |>1-\frac{1}{k}\} )\subset
\overline{W_0}\}.
\]
Each $E_k$ is a closed hence measurable set and $E=\bigcup_{k\geq
1}E_k$. By countable additivity there is a set $E_k$ that has
positive measure. Set $E_0=E_k$, $\alpha=\frac{1}{k}$,
$r=1-\frac{1}{k}$. Then  $u (\Gamma_r(E_0,\alpha ))\subset
\overline{W_0}$, and that establishes the Step 1.

\medskip\noindent{\slsf Step 2.} {\sl If for $e^{i\theta}$ in a set $E_0$ of positive length
in $\d\Delta$ one has that $u^\ast (e^{i\theta})\in V_{j}\cap W$ for some relatively compact
$W\subset X\setminus\bigcup_{l=0}^{j-1}V_l$ then $u(U)\subset V_{j}$ for some non-empty
open subset $U\subset\Delta$.}

\smallskip
One can replace $E_0$ by a subset small enough if needed, so that $\Gamma_r(E_0,\alpha )$
is a simply connected domain with rectifiable boundary $\gamma$. This is well explained
in \S 2 of Chapter X in \cite{Go}.
Let $\chi$ be a conformal mapping from the unit disc $\Delta$ onto $\Gamma_r(E_0,\alpha )$.
This map extends continuously to the closed unit disc.
Set $F=\{ e^{i\nu}\in \d\Delta : \chi (e^{i\nu})\in E_0\}$. Then $F$ is a set of positive length
in $\d\Delta$ by F. and M. Riesz's theorem on conformal mappings, see Theorem VIII.26 p. 318
in \cite{T}.

For $\epsilon >0$ small enough, define $\rho_\epsilon$ on $\Delta$
by: $\rho_\epsilon (\zeta) = \rho\circ u\big( (1-\epsilon )\chi
(\zeta )\big)$. Note that $(1-\epsilon )\Gamma_r(E_0,\alpha )\subset
\Gamma_{(1-\epsilon )r}(E_0,\alpha )$. Then, the functions
$\rho_\epsilon$ are subharmonic functions on $\Delta$ that are
uniformly bounded from above, and for any $e^{i\nu}\in F$, by upper semi-continuity of
$\rho$,
$\rho_\epsilon (e^{i\nu})$ tends to $-\infty$ as $\epsilon\to 0$.
Indeed $(1-\epsilon )\chi (e^{i\nu})$ approaches $\chi (e^{i\nu})\in
E_0$ radially. So if $e^{i\nu}\in F$, $\chi
(e^{i\nu})=e^{i\theta}\in E_0$, then $u((1-\epsilon )\chi
(e^{i\nu}))$ tends to $u^\ast (e^{i\theta})$ as $\epsilon \to 0$. By
the mean value property, $\rho_{\eps}$ tends to $-\infty$ uniformly
on compact subsets of $\Delta$. Since $\rho\circ u (\zeta) =
\rho_\epsilon (\chi^{-1}(\frac{\zeta}{1-\epsilon}))$ (for $\eps <
1-|\zeta |$), $\rho\circ u  \equiv -\infty$ on $\Gamma_r(E_0,\alpha
)$, so $u(\Gamma_r(E_0,\alpha ))\subset V_{j}\cap W$.

\medskip\noindent{\slsf Step 3.} {\sl In the conditions of Step 2 there exists a connected dense
open subset
$A\subset \Delta$ such that $u(A)\subset V_j$.}

\smallskip\noindent
Let $A$ to be the set of points $\zeta\in\Delta$ such that $u$ maps
a neighborhood of $\zeta$ into $V_j$, let $\d A$ be its boundary.
Clearly $A$ is an open subset of $\Delta$. For each $l\leq j$, let
$S_l=\{\zeta\in \d A ;~u(\zeta )\in V_l\}$. If $\zeta\in S_l$, there
is a $J$-psh function $\rho$ defined near $u(\zeta )$, that near
this point is $-\infty$ exactly on $V_l$. Then $\rho\circ u=-\infty$
near $\zeta$ on $S_l$, but is not identically $-\infty$ since
$\zeta\in \d A$ (more precisely: if $l<j$ because $\zeta$ is in the
closure of $A$, if $l=j$ because $\zeta\notin A$). Hence each $S_l$
is a polar subset of $\Delta$ and so is their union that is $\d A$.
Since $\d A$ is a closed polar subset of $\Delta$, $A$ is  a
connected dense open subset of $\Delta$, as claimed.

\medskip By continuity, $u(\Delta )\subset
\bigcup_{l=o}^jV_l\subset V$. The proof of Theorem \ref{geo-incl} is
therefore completed.

\smallskip\qed

\medskip\noindent{\slsf 3.3. The diagonal and proof of Privalov-type
Theorem}

\smallskip
Let us derive from Theorem \ref{geo-incl} the Privalov-type Theorem
\ref{priva} from the Introduction. Let $(\Omega_1,\rho_1)$ and
$(\Omega_2,\rho_2)$ be relatively compact domains containing
$\overline{u_1(\Delta)}$ and $\overline{u_2(\Delta)}$ together with
strictly $J$-psh functions, which clearly could be supposed to be
bounded from above. Then $\Omega\deff(\Omega_1\times\Omega_2, \tilde
J\deff J\oplus J, \lambda_1\deff\rho_1 + \rho_2)$ will be the same
type data for $u\deff (u_1,u_2):\Delta\to (X\times X,\tilde J)$, \ie
$u$ is also bounded. Remark further that $u:= (u_1,u_2): \Delta\to
\Omega\subset X\times X$ is $\tilde J$-holomorphic.

\smallskip With $\Omega $ and $\rho$ just defined we consider the intersection
$V\deff \dd\cap \Omega $ of the diagonal $\dd$ of $X\times X$ with
$\Omega$. This is a $\tilde J$-complex submanifold of $\Omega$. Note
that if for at least one $e^{i\theta}$ the boundary values
$u_1^*(e^{i\theta})$ and $u_2^*(e^{i\theta})$ do coincide then $V$
is nonempty.

\smallskip
By  the theorem of Elkhadhra  $V$ is locally complete polar and
therefore  our problem is reduced to the study of boundary values of
a $J$-holomorphic map lying on a locally complete pluripolar subset
of an almost complex manifold. Theorem \ref{geo-incl} is now
applicable to $u$ and gives us that $u(\Delta)\subset V$. Theorem
\ref{priva} follows.

\section[4]{Blaschke Condition for a Single Pseudoholomorphic Mapping}

\smallskip\noindent{\slsf 4.1. The Blaschke Condition: proof of Theorem \ref{bla}.}

\smallskip Let $\lambda_0$ be a strictly $J$-psh
function defined on a neighborhood of $\overline {u(\Delta )}$.
Using Chirka's function, we claim that there exists a $J$-psh
function $\lambda$ defined on $\Omega$, bounded from above and such
that near $p$, using local coordinates: $\lambda (z)\leq \log|z-p|$.
Indeed, using local coordinates near $p$ such that $J(p)=J_{st}$,
set
$$\lambda (z)= \chi (z)(\log |z-p| + A|z-p|)+B\lambda_0-C~,$$
where $\chi$ is a cut off function equal to 1 on a neighborhood of
$p$. $A$ is first taken large enough in order that
$\log|z-p|+A|z-p|$ is $J$-psh near $p$, then $B$, then $C$ have just
to be taken large enough. We set $\rho = \lambda \circ u$. The
Blaschke condition follows from the $\log$ singularity due to the
following two facts:

\medskip\noindent
(a) If $\rho$ is a subharmonic function defined near $\zeta_0$ in ${\bf C}$
and if for $\zeta$ close to $\zeta_0$ one has
$\rho (\zeta )\leq \log |\zeta -\zeta_0 |~-~C$ for some constant $C$,
then $\Delta (\rho )(\{\zeta_0 \}) \geq 2\pi$, where $\Delta (\rho )$ denotes the
positive measure that is the distributional Laplacian of $\rho$.
An immediate proof for this classical fact consists in using that,
near $p$, $\rho$ is the decreasing limit as $n\to \infty$ of the functions
$\rho_n={\rm max}~\big(\rho, (1-\frac{1}{n}
)\log |\zeta -\zeta_0|~-n\big)$, for which
obviously $\Delta (\rho _n) (\{ \zeta_0\}) = 2\pi (1-\frac{1}{n})$.

\medskip\noindent
(b) If $\rho$ is a subharmonic function on $\Delta$ that is bounded
from above, then
\begin{equation}
\label{int-bound}
\int_{\Delta} (1-|\zeta |)\Delta \rho (\zeta )~d\xi d\eta
<+\infty~.
\end{equation}
This also is classical and it follows from the representation formula
\begin{equation}
\eqqno(int-ineq)
\rho (0) = \frac{1}{2\pi }\int_0^{2\pi}\rho (e^{i\theta}) d\theta
-\frac{1}{2\pi} \int_{\Delta} \Delta (\rho )(\zeta )~\log\frac{1}{|\zeta|}d\xi d\eta ,
\end{equation}
(see e.g. \cite{Ga} XV.2, page 294), and from the fact that $\log\frac{1}{|\zeta|}\sim 1 - |\zeta|$ near the
boundary of $\Delta$.

\medskip We can now conclude the proof. We apply (a) and (b) with
$\rho = \lambda \circ u$. For any $\zeta_0$
in $\Delta$ such that $u(\zeta_0)=p$, we have for $\zeta$ close to
$\zeta_0$ that $|u(\zeta)-p|\leq C|\zeta -\zeta_0|$. Hence, $\rho\leq
\log |\zeta -\zeta_0|+ C$, therefore, by (a), $\Delta (\rho )(\{ \zeta_0 \})\geq 2\pi$.
Then (\ref{int-bound}) gives us the desired Blaschke condition.

\smallskip\qed

\smallskip\noindent{\slsf 4.2. An example of a ``bad diagonal''}

\smallskip In \S 4.1, it has been essential that the Chirka function has a
logarithmic singularity, not only a pole. In \cite{Ro2}, it has been
shown that $J$-holomorphic curves are not in general $-\infty$ of
$J$-psh functions with logarithmic singularity. In \S 3.3 we applied
Elkhadhra's Theorem to a very special case,
the diagonal in the
product of the space with itself. Unfortunately, even in that
special case, there may not exist plurisubharmonic functions with
logarithmic singularity. We are thus unable to extend Theorem
\ref{bla} to establish that the set of points where two
$J$-holomorphic maps coincide, is a Blaschke sequence, under
appropriate boundedness assumption.

\medskip We use the example given in \cite{Ro2}. In that
example a smooth almost complex structure $J$ on $\cc_{z_1} \times
\cc_{z_2}$ is constructed such that:

\smallskip

(i) $\{ z_2=0\}$ is $J$-holomorphic, but is not the $-\infty$ of a
$J$-psh function with logarithmic singularity, not even locally.

\smallskip

(ii) For every $J$-holomorphic disc $\zeta \mapsto  (u_1(\zeta ),u_2(\zeta ))$, $u_1$ is holomorphic.

\smallskip

(iii) $\zeta\mapsto (u_1(\zeta ),0)$ is $J$-holomorphic if $u_1$ is holomorphic. (But $\zeta\mapsto (u_1(\zeta ),a)$,
for fixed $a\not=0$ is not).

\smallskip\noindent
We now consider the product structure $J\oplus J$ on $\cc^2\times \cc^2$, and we let $\dd$ be the diagonal.

\begin{prop}
There is no $J\oplus J$-plurisubharmonic  function with logarithmic singularity along $\dd$.

\end{prop}

\noindent Proof. The imbedding from $(\cc^2,J)$ into $(\cc^2\times
\cc^2, J\oplus J)$ defined by
\[
\phi : (z_1,z_2)\mapsto  (z_1,z_2,z_1,0)
\]
is $(J, J\oplus J)$-holomorphic.
Indeed, for any $J$-holomorphic map $u$, $\phi\circ u$ is $J\oplus J$-holomorphic, as it follows
immediately from (ii) and (iii).
\smallskip\noindent
If $\rho$ is a $J\oplus J$-plurisubharmonic function on
$(\cc^2\times \cc^2,J\oplus J)$, that we shall suppose to be defined
near a point $(z_1,0,z_1,0)$, then $\rho\circ \phi$ is $J$-psh. If
$\rho$ had a logarithmic singularity along $\dd$ then $\rho\circ
\phi$ would have a logarithmic singularity along $z_2=0$, which is
not possible.

\section[5]{Appendix}

\noindent{\slsf 5.1. Existence of $J$-Complex Curves in Continuous Structures}

We consider a continuous operator valued function $J$ in the
unit ball $B$ of $\rr^{2n}$. I.e., $J:B\to ${\slsf
Mat}$(2n\times 2n,\rr)$. $J$ is an almost complex structure if
$J^2(x)\equiv -\id $. We will use the notation $W^{1,p}$ for the Sobolev spaces
of $L^p$ functions with first derivatives in $L^p$.

\begin{defi}
A $W^{1,1}_{loc}\cap\calc^0$-map $u:\Delta\to B$ is said to be
$J$-holomorphic if it a.e. satisfies
\begin{equation}
\eqqno(J-hol1)
\frac{\d u}{\d x} + J(u(z))\frac{\d u}{\d y} =0.
\end{equation}
\end{defi}
The image $C=u(\Delta)$ is called a $J$-complex disk.
$J_u(z):=J(u(z))$ can be considered as a matrix valued function on
the unit disk and therefore can be viewed as a complex linear
structure on the trivial bundle $E:=\Delta\times\rr^{2n}$. I.e.,
$J_u(z)\in \calc^0(\Delta , End(\rr^{2n}))$ satisfying $J(z)^2\equiv
-\id $. The mapping $u$ can be viewed as a section of this bundle. It is
known that if $J\in\calc^0$ then $J$-holomorphic maps belong to
$W^{1,p}_{loc}$ for all $p<\infty$ (see \cite{ISh1} Lemma 1.2.2) and
therefore to $\calc^{\alpha}$ for all $\alpha < 1$ by Sobolev
imbedding theorem. $J$-complex curves in continuous structures have
some nice properties. For example, a Gromov compactness theorem is
valid for them, see \cite{ISh2}. Hovewer the question of existence
of complex curves in continuous structures to our knowledge was
never discussed in the literature, therefore we do this here.

\smallskip Equation \eqqref(J-hol1) can be rewritten as

\begin{equation}
\eqqno(J-hol2)
\frac{\d u}{\d\bar z} + \bar Q_J(u(z))\frac{\d u}{\d z} = 0,
\end{equation}
where $\frac{\d u}{\d\bar z}=\frac{1}{2}(\frac{\d u}{\d x}+J\st
\frac{\d u}{\d y})$, $\frac{\d u}{\d z}=\frac{1}{2}(\frac{\d u}{\d
x}-J\st \frac{\d u}{\d y})$ and

\begin{equation}
\eqqno(matrix-q)
\bar Q(J(z))=[J(z)+J\st ]^{-1}[J(z)-J\st].
\end{equation}
Remark that $\bar Q$ anti commutes with $J\st$ and therefore is a
$\cc$-antilinear operator. Thus \eqqref(J-hol2) can be understood as
an equation for $\cc^n$-valued map (or section) $u$. Usually it is
better to consider the conjugate operator $Q$ and write
\eqqref(J-hol2) in the matrix form

\begin{equation}
\frac{\d u}{\d\bar z} + Q_J(u(z))\overline{\frac{\d u}{\d z}} = 0.
\eqqno(conjug-q)
\end{equation}
We shall look for the solutions of \eqqref(conjug-q) in the form

\begin{equation}
u(z) = - T_{CG}Q(J_u(z))\overline{\frac{\d u}{\d z}}~+H =:\Phi(u)~+H,
\eqqno(phi)
\end{equation}
where $T_{CG}$ is the Cauchy-Green operator (convolution with
$\frac{1}{\pi z}$), and $H$ is a holomorphic function. We shall scale
the norm on $W^{1,p}(\Delta ,\cc^n)$, so that for any $f\in
W^{1,p}(\Delta ,\cc^n)$
\begin{equation}
\eqqno(scale)
\norm{f}_{L^{\infty}(\Delta)}\leq \|
f\|_{W^{1,p}(\Delta)}.
\end{equation}
\begin{lem}
\label{fixed}
Suppose that $J(0)=J\st $, so $Q(J(0))=0$. Let $2<p<\infty$, there exists $q>0$ such that
if $\norm{Q(J)}_{L^{\infty}(B)}=q<<1$, then
for $a,b\in\cc^n$ small enough, the operator $\Phi_{a,b}$ defined by
\begin{equation}
\Phi_{a,b}(u) = \Phi(u) - \Phi(u)[0] -
2z\left(\Phi(u)[1/2]-\Phi(u)[0]\right) + a + 2z(b-a)
\eqqno(phi-ab)
\end{equation}
is an  operator from the closed unit ball $\calb\subset
W^{1,p}(\Delta ,\cc^{n})$ into itself, that has a fixed point.
\end{lem}
\smallskip\noindent{\sl Proof.} Unfortunately the lemma is not obtained directly by
using the contraction principle, instead one uses an argument
already in \cite{NW}. The proof consists in two steps:

\smallskip\noindent
(i) Let $\calb_\infty$ be the closed unit ball in $L^\infty(\Delta
,\cc^n)$. Note that $\calb \subset \calb_\infty$ due to
\eqqref(scale). For any $v\in\calb_\infty$, one defines the operator
\begin{equation}
\eqqno(phi-abv) \Phi_{a,b}^v(u)= \Phi^v(u) - \Phi^v(u)[0] -
2z\left(\Phi^v(u)[1/2]-\Phi^v(u)[0]\right) + a + 2z(b-a) ~,
\end{equation}
where
\begin{equation}
\eqqno(phi-v)
\Phi^v(u) \deff - T_{CG}Q(J_v(z))\overline{\frac{\d
u}{\d z}},
\end{equation}
thus taking $J(v(z))$ instead of $J(u(z))$. By the contraction
principle, one proves that this operator has a unique fixed point
$u=T(v)$ in $\calb$.

\smallskip\noindent
(ii) Then one proves that $T$ restricted to $\calb$, itself has a
fixed point $u$, by applying Schauder fixed point theorem.
Then $\Phi_{a,b}^u(u)= u$, i.e. $\Phi_{a,b} (u)=u$,  and that establishes the Lemma.

\bigskip\noindent
(i) Let $C_{p}$ be the norm  of the Cauchy-Green operator as a
linear map from $L^{p}(\Delta)\to W^{1,p}(\Delta)$. If
$q<\frac{1}{8C_{p}}$ then for $u_1,u_2\in \calb$
\[
\norm{\Phi^v_{a,b} (u_1)-\Phi^v_{a,b}(u_2)}_{W^{1,p}}\le
4\norm{T_{CG}Q(J_{v})\left[\overline{\frac{\d u_1}{\d
z}}-\overline{\frac{\d u_2}{\d z}}\right]}_{W^{1,p}} \leq
\]
\[
\leq 4qC_{p}\norm{u_1-u_2}_{W^{1,p}}\leq
\frac{1}{2}\norm{u_1-u_2}_{W^{1,p}}.
\]
Since $\Phi^v_{a,b}(0)=a+2z(b-a)$, for $a$ and $b$ small enough
$\Phi^v_{a,b}$ is a contracting map from $\calb$ into itself. Hence it has a unique fixed point
$u=T(v)$.

\bigskip\noindent
(ii) Let us prove now that $v\mapsto u=T(v)$ defines a continuous
(in $L^{\infty}$-topology) map from $\calb_\infty$ to
$\calb\subset\calb_{\infty}$. Indeed, we have
\[
u=\Phi_{a,b}^{v}(u)= \Phi(u) - \Phi(u)[0] -
2z\left(\Phi(u)[1/2]-\Phi(u)[0]\right) + a + 2z(b-a) ~ \]
and thus
\[
\norm{u}_{W^{1,p}} = \norm{\Phi_{a,b}^{v}u}_{W^{1,p}} \le
4qC_p\norm{u}_{W^{1,p}} + \norm{a + 2z(b-a)}_{W^{1,p}}.
\]
We get therefore
\begin{equation}
\eqqno(bound)
\norm{u}_{W^{1,p}}\le \frac{1}{1-4qC_p}\norm{a +
2z(b-a)}_{W^{1,p}},
\end{equation}
which means that for $a$ and $b$ small enough $T$ maps
$\calb_{\infty}$ to $\calb$.

\smallskip
Let $(v_n)$ be a converging sequence in $\calb_{\infty}$ with limit
$v$, and $u_n=Tv_n$. \eqqref(bound) tells us that $u_n$ are bounded
in $W^{1,p}(\Delta)$. Therefore the sequence ${\frac{\d u_n}{\d z}}$
is a bounded in $L^p(\Delta)$. Let $u=Tv$, \ie $u$ is the unique
fixed point of $\Phi^v_{a,b}$. Write:
\[
||u-u_n||_{W^{1,p}}=||\Phi_{a,b}^{v}(u)-\Phi_{a,b}^{v_n}(u_n)||_{W^{1,p}}\le
4||\Phi^{v}(u) - \Phi^{v_n}(u_n)||_{W^{1,p}} =
\]
\[
= 4\norm{T_{CG}Q(J_{v})\left[\overline{\frac{\d u}{\d
z}}-\overline{\frac{\d u_n}{\d z}}\right]}_{W^{1,p}} +
4\norm{T_{CG}\left[Q(J_{v})-Q(J_{v_n})\right]\overline{\frac{\d
u_n}{\d z}}}_{W^{1,p}} \le
\]
\[
\le 4qC_p\norm{u-u_n}_{W^{1,p}} +
4\norm{T_{CG}\left[Q(J_{v})-Q(J_{v_n})\right]\overline{\frac{\d
u_n}{\d z}}}_{W^{1,p}}.
\]
Therefore
\begin{equation}
||u-u_n||_{W^{1,p}} \le
\frac{4}{1-4qC_p}\norm{T_{CG}\left[Q(J_{v})-Q(J_{v_n})\right]\overline{\frac{\d
u_n}{\d z}}}_{W^{1,p}} \to 0,
\end{equation}
because $Q(J_{v_n})$ converges uniformly to $Q(J_v)$. The continuity
of $T$ is proved.

\smallskip
Take a closure $\bar\calb$ of $\calb$ in $\calb_{\infty}$. This is a
convex compact subset of $\calb_{\infty}$ in $L^{\infty}$-topology.
By the Schauder fixed point Theorem applied to $T|_{\bar\calb}$, $T$
has a fixed point $u\in\bar\calb$ and this $u$, as we had proved,
belongs to $\calb$ in fact.

\smallskip\qed

\begin{corol}
Let $J$ be a continuous almost complex structure in $\rr^{2n}$. Then
for any sufficiently close points $a$ and $b$ there exists a
$J$-holomorphic map $u:\Delta\to\rr^{2n}$ passing through them.
\end{corol}
The result follows from Theorem \ref{fixed}.
Indeed, for the fixed point $u$ of $\Phi_{a,b}$
we get
$\frac{\partial u}{\partial \overline z}+Q_J(u)\overline{\frac{\partial u}{\partial z}}=0$,
$u(0)=a$ and $u(\frac{1}{2})=b$. The needed smallness of $Q$ (equivalently
smallness of $J-J_{st}$), given $J(0)=J_{st}$, is obtained by simple rescaling.

\smallskip\qed

\begin{rema}\rm
The standard maximum principle for $J$-psh functions
in case of continuous $J$'s follows immediately.
\end{rema}

\smallskip\noindent{\slsf 5.2. Plurisubharmonicity and Pluripolarity in Lipschitz structures.}

\smallskip\noindent If an almost complex structure $J$ is not differentiable, the operator
$dd^c_J$ (whose definition requires one differentiation of $J$) does
not have an immediate meaning. For studying Lipschitz structures,
instead of trying to define $dd^c_J$ we shall proceed by
approximation.

\begin{lem}
\label{lipJ} Let $J$ be a Lipschitz-continuous almost complex
structure defined on a (smooth) manifold $X$. Let $J_k$ be a
sequence of $\calc^1$ almost complex structures on $X$, with
$\calc^1$-norms uniformly bounded, and converging uniformly to $J$
as $k\to\infty$. Let $\rho$ be a $\calc^2$ function defined on $X$.
If $\rho$ is $J_k$-psh, for $k$ large enough, then $\rho$ is
$J$-psh.
\end{lem}

\smallskip\noindent {\sl Proof.} The problem is local, we can always
restrict our attention to a relatively compact domain in $X$, and we
equip $X$ with some Riemannian metric. We shall use the following
characterization of $\calc^1$ (but not $\calc^2$) subharmonic
functions $\mu$, on the unit disc $\Delta$:
\begin{equation}
\eqqno(c1-char)
-\int_{\Delta}d\varphi \wedge d^c\mu~~\geq 0~,
\end{equation}
for all $\varphi\in \calc^1_0(\Delta )$, with $\varphi \geq 0$.
(Here $d^c=d^c_{J_{st}}$). The above condition gives $dd^c\mu \geq
0$, in the sense of distributions.

\bigskip
Let $u:\Delta \to (\rr^{2n},J)$ be a $J$-holomorphic map from
$\Delta$ into $B$, so $u$ is of class $\calc^{1,\beta}$ for all
$\beta <1$. We can assume $u$ to be $\calc^{1,\beta}$ up to the
boundary. We wish to show that $\rho\circ u$ is subharmonic. Since
$u$ is $J$-holomorphic $d^c\rho\circ u=u^\ast d^c_J\rho$ (by the
simple rules of differentiation and commutation of $d$ with the
action of the almost complex structures). So, we need to show that
for any non-negative $\calc^1$ function $\varphi$ with compact
support in $\Delta$,
\begin{equation}
\label{lipJ0}
-\int_{\Delta}d\varphi \wedge u^\ast d^c_J\rho~~\geq 0~.
\end{equation}

\bigskip\noindent
By hypothesis we have $dd^c_{J_k}\rho (v,J_kv)\geq 0$, for all tangent vector $v$ at any point
$p$. Let $\epsilon_k= {\rm Sup}_p|J_k(p)-J(p)|$, so $\epsilon_k\to 0$ as $k\to \infty$.
Since $\rho$ is of class $\calc^2$ and since the $\calc^1$ norms of $J_k$ are uniformly bounded,
for some constant $K$ we have

\begin{equation}
\label{lipJ1}
dd^c_{J_k}\rho (v,Jv)\geq -K\epsilon_k \| v\|^2.
\end{equation}
\smallskip\noindent
The inequality (\ref{lipJ1}) has a clear geometric meaning: there
exists a constant $K_1>0$ such that for any germ of embedded
$J$-holomorphic curve $u$ in $X$, oriented by $J$, $dd^c_{J_k}$
induces on that curve a measure $\lambda_k$ that satisfies
\begin{equation}
\label{lipJ2}
\lambda_k\geq - K_1\epsilon_k dm~,
\end{equation}
where $dm$ denotes Euclidean area (Hausdorff) measure on the curve.

\smallskip
We now prove (\ref{lipJ0}). By uniform convergence of $J_k$,
$$ -\int_{\Delta}d\varphi \wedge u^\ast d^c_J\rho ~=~
{\rm lim}_{k\to \infty}  -\int_{\Delta}d\varphi \wedge u^\ast d^c_{J_k}\rho~.$$
Set $I_k=-\int_{\Delta}d\varphi \wedge u^\ast d^c_{J_k}\rho$. Because
$u$ is not $\calc^2$, the differential form $ u^\ast d^c_{J_k}\rho$
may fail to be $\calc^1$, since pull back of a differential form by a map uses a derivative of the map.
However Stokes formula can still be applied and it gives us
$I_k=\int_{\Delta}\varphi  u^\ast dd^c_{J_k}\rho$, where the right hand side clearly makes sense.
See the preliminary remark in the proof of Lemma 1.2 in \cite{IR}.
Then (\ref{lipJ2}) yields $I_k\geq -(K_1~{\rm Sup}~\varphi )~M\epsilon_k~,$
where $M$ is the area of $u(\Delta )$. Therefore (\ref{lipJ0}) follows and the Lemma is proved.

\bigskip\noindent
We can now generalize the plurisubharmonicity result for the Chirka function to the
case of Lipschitz structures.
\bigskip\noindent
\begin{lem}
\label{chiL} Let $J$ be a Lipschitz continuous almost complex
structure defined near 0 in $\cc^n$. If $J(0)=J_{st}$, then for
$A>0$ large enough, the function $\rho$ defined by $\rho (z)=\log
|z|~+~A|z|$ is $J$-psh.
\end{lem}

\smallskip\noindent {\sl Proof.}
We can approximate $J$ uniformly by a sequence $(J_k)$ of almost complex structures
with $J_k(0)=J_{st}$, and with uniformly bounded $\calc^1$ norms.
By Lemma 1.4 in \cite{IR}, for $A$ large enough, for $k$ large enough
$\rho$ is $J_k$-psh (that $A$ can be chosen independently of $k$ (large) is clear from the last lines
of the proof of the lemma in \cite{IR}). Since plurisubharmonicity needs to be proven only
off the $-\infty$ set, and since $\rho$ is smooth except at 0, we can apply Lemma \ref{lipJ}.
That establishes Lemma \ref{chiL}.

\smallskip\qed

\ifx\undefined\bysame
\newcommand{\bysame}{\leavevmode\hbox to3em{\hrulefill}\,}
\fi

\def\entry#1#2#3#4\par{\bibitem[#1]{#1}
{\textsc{#2 }}{\sl{#3} }#4\par\vskip2pt}

\end{document}